\documentclass[letterpaper, 10 pt, conference]{ieeeconf}  % Comment this line out if you need a4paper

\IEEEoverridecommandlockouts                              % This command is only needed if 
\usepackage{Preamble}
\title{\LARGE \bf
	Risk-aware stochastic control of a sailboat 
}

\author{MingYi Wang$^{1},$ Natasha Patnaik$^{2},$ Anne Somalwar$^{3},$ Jingyi Wu$^{4},$ and Alexander Vladimirsky$^{5}$% <-this % stops a space
	\thanks{*This work was supported in part by NSF Division of Mathematical Sciences (DMS) under Awards 1645643 and 2111522 as well as 
	by the Air Force Office of Scientific Research under Award FA9550-22-1-0528.}% <-this % stops a space
	\thanks{$^{1}$MingYi Wang is with the Center for Applied Mathematics,
		Cornell University, Ithaca, NY 14853, USA
		{\tt\small mw929@cornell.edu}}%
	\thanks{$^{2}$Natasha Patnaik is with the Department of Computational Applied Mathematics and Operations Research, Rice University,
		Houston, TX 77005, USA
		{\tt\small np37@rice.edu}}%
	\thanks{$^{3}$Anne Somalwar is with the graduate group in Applied Mathematics and Computational Science, University of Pennsylvania,
		Philadelphia, PA 19104, USA
		{\tt\small somalwar@sas.upenn.edu}}%
	\thanks{$^{4}$Jingyi Wu is with the Center for Data Science, New York University,
		New York City, NY 10012, USA
		{\tt\small  jw8535@nyu.edu}}%
	\thanks{$^{5}$Alexander Vladimirsky is with the Department of Mathematics,
		Cornell University, Ithaca, NY 14853, USA
		{\tt\small vladimirsky@cornell.edu}}%
}

\begin{document}
	\maketitle
	\thispagestyle{empty}
	\pagestyle{empty}
	
%%%%%%%%%%%%%%%%%%%%%%%%%%%%%%%%%%%%%%%%%%%%%%%%%%%%%%%%%%%%%%%%%%%%%%%%%%%%%%%%
\begin{abstract}
Sailboat path-planning is a natural hybrid control problem (due to continuous steering and occasional “tack-switching” maneuvers), with the actual path-to-target greatly affected by stochastically evolving wind conditions.  Previous studies have focused on finding risk-neutral policies that minimize the expected time of arrival.
In contrast, we present a robust control approach, which maximizes the probability of arriving before a specified deadline/threshold.
Our numerical method % COMMENT REMOVED 
recovers the optimal risk-aware (and threshold-specific) policies for all initial sailboat positions and a broad range of thresholds simultaneously.  This is accomplished by solving two quasi-variational inequalities based on second-order Hamilton-Jacobi-Bellman (HJB) PDEs with degenerate parabolicity.  
% COMMENT REMOVED
Monte-Carlo simulations show that risk-awareness in sailing is particularly useful when a carefully calculated bet on the evolving wind direction might yield a reduction in the number of tack-switches.
% COMMENT REMOVED
\end{abstract}

% COMMENT REMOVED
%===============================================================================

% !TEX root = Main_text.tex

\section{Introduction}\label{sec:Intro}
Sail-boat racing is one of the many areas where game-theoretic and control-theoretic tools are valuable in improving the competitive performance. 
The uncertainty in weather patterns gives rise to hybrid stochastic control models with many 
reasonable choices of performance measures to optimize.
% COMMENT REMOVED
The previously developed methods have focused on risk-neutral optimization (e.g., minimizing the expected time to destination) \cite{Ferretti2019, cacace_stochastic_2019, Miles2021}.  In contrast, here we focus on  maximizing the probability of desirable outcomes (e.g., arriving prior to a specified deadline).  Our {\em risk-aware} approach addresses a notion of robustness very different from the traditional $H^\infty$ control \cite{bacsar2008h} and has important advantages for many applications.  
% COMMENT REMOVED
Indeed, it has been already successfully used in piecewise-deterministic Markov processes \cite{CarteeVlad_UQ} and in bang-bang stochastic control models of adaptive drug therapies \cite{WangScottVlad_2023}.  
Unlike the typical {\em risk-averse} approaches \cite{WangChapman_2022},
the method that we develop here for the hybrid control setting allows finding the optimal control policies for a large set of starting sailboat positions and 
{\em a range of deadlines} simultaneously.  This is accomplished in the general framework of dynamic programming and requires solving a pair of quasi-variational inequalities on a 
% COMMENT REMOVED
3D computational domain.

The hybrid nature of this problem is due to ``tacking'': to travel upwind, sailors must use a zigzag pattern, periodically swinging the bow to the other side of the wind.
Each such ``tack-switch'' incurs a significant slow down while the wind pushes against the boat.
We adopt a commonly used simplified model which assumes
% COMMENT REMOVED
that the boat's velocity vector can be changed instantaneously (choosing among all directions available in its current tack) but a switch to the opposite tack incurs a fixed time-penalty.  

Optimization of sailboat routing 
% COMMENT REMOVED
is
a topic of increasing mathematical interest. 
 In \cite{philpott_simulation_2004} the task of minimizing the expected time to target was considered in a discrete setting, with a discrete-time Markov chain modeling occasional changes in weather conditions.  The idea was extended to continuous-time Markov chains in \cite{vinckenbosch_stochastic_2012}, with a tack-switching curve defined in the state space to encode the optimal policy.  
% COMMENT REMOVED
In \cite{Ferretti2019}, this switching curve was found using dynamic programming for indefinite-horizon hybrid control problems, 
but under the assumption that the wind direction stays constant for the duration of each tack switch.
In \cite{Miles2021}, it was shown how this assumption can be avoided, yielding an improvement in control policies.

We start by introducing the % COMMENT REMOVED
hybrid dynamics in Section \ref{sec:Sys_Dynamics}, 
and describe both the risk-neutral and risk-aware optimal control problems in Section \ref{sec:Sto_OC}.
Our numerical approach to the latter is presented in Section \ref{sec:Num_Implementaion}, 
followed by the summary of computational experiments in Section \ref{sec:Num_Experiments}.
We conclude by listing directions for future work in Section \ref{sec:Conclusion}.

% !TEX root = Main_text.tex

\section{System Dynamics}\label{sec:Sys_Dynamics}
\begin{figure*}[ht]
	\centering
	\subfigure[The polar speed plot $f(u)$ used in this work (the same as Fig.~1(a) in \cite{Miles2021}). Note that $f=0$ at $u=0^\circ$.]{\label{fig:speed_and_tacks_a}\includegraphics[scale=0.38]{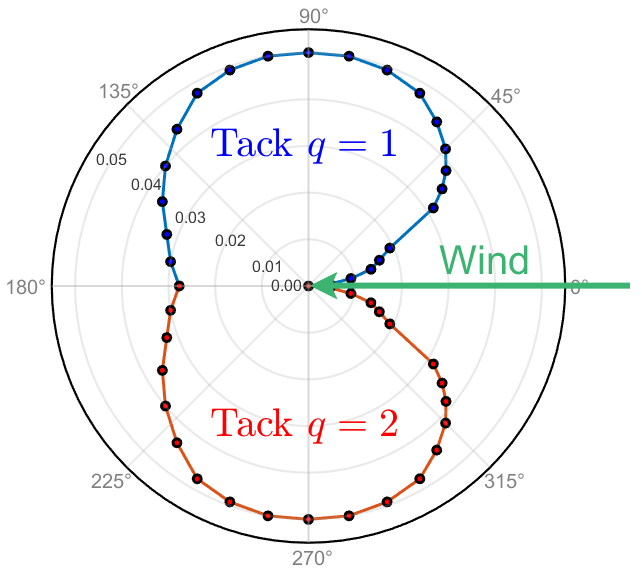}}
	\subfigure[System setups in the polar coordinate centered at $\tgt$ for different tacks.]{\label{fig:speed_and_tacks_b}
	   \centering
		% Tack 1
		 \begin{tikzpicture}[scale=0.9]
			
			 % Drawing half circle
			  \coordinate (A) at (4.553, 1.224);
			  \coordinate (B) at (5.447,0.776);
			
			  \coordinate (C) at ($(A)!.5!(B)$);

			   \draw (A) -- (B);
			   \pgfmathanglebetweenpoints{\pgfpointanchor{C}{center}}{\pgfpointanchor{A}{center}}
			   \let\StartAngle\pgfmathresult
			   \pgfmathanglebetweenpoints{\pgfpointanchor{C}{center}}{\pgfpointanchor{B}{center}}
			   \let\EndAngle\pgfmathresult

			   \draw[fill = blue!40] (A) let \p1 = ($(C) - (A)$),
			   \n1 = {veclen(\x1, \y1)}
			   in 
			   arc [start angle= \StartAngle , end angle= \EndAngle , radius=\n1];
			
			 \draw[dashed,line width = 0.3 mm] (3.5,5.5) -- (5,1);
			 \draw (3.5,5.5) node[draw,circle,fill=blue!20,inner sep=5pt]{$\tgt$};
			 \node[] at (1.9,6.3) {\large \textcolor{blue}{\textbf{Tack $q=1$}}};
			
			 \draw [decorate, line width=0.5mm,
			     decoration = {calligraphic brace, raise = 20pt}] (3.5,5.5) -- (5,1)
			     node[pos=0.45, right=30pt,black]{$r$};
			
			 %Defining points
			 %Point A is the location of the boat
			 %Point B is the tip of the dr vector
			 %Point C is the tip of the v vector
			 %Point D is the foot of the wind vector
			 %Point E is the tip of the wind vector
			  \tkzDefPoints{5/1/A,4.44991/2.65023/B,1/1.5/C,1/3/D,6/0.5/E}
			
			 %Drawing vectors / labeling
			  \tkzDrawSegments[vector style, violet!50,line width = 0.5 mm](A,B)
			  \tkzDrawSegments[vector style,violet,line width = 0.5 mm](A,C)
			  \tkzDrawSegment[vector style,violet!50,line width = 0.5 mm](B,C)
   			  \tkzDrawSegment[vector style, green!40!gray,line width = 0.5 mm](D,E)
			  \tkzMarkAngle[->, blue, size = 1.5,line width = 0.4 mm](D,A,C)
			  \tkzLabelAngle[blue, pos=1.2](D,A,C){$u$}
			  \tkzMarkAngle[->, green!40!gray, size = 2.7,line width = 0.4 mm](B,A,D)
			  \tkzLabelAngle[green!40!gray, pos=2.5, left = 16pt](B,A,D){$\theta$}
			  \tkzLabelSegment[violet!40, right = 1pt](A,B){$dr$}
			  \tkzLabelSegment[violet,swap, pos = 0.92, below left=  0.1 pt](A,C){$\vec{v}$}
			
			  \tkzLabelSegment[violet!50, pos = 0.15, above = 2pt](C,B){$d\theta$}
			
			  \tkzLabelSegment[green!40!gray, pos = 0, above right= 0.5 pt ](D,E){\textbf{Wind}}
			  
			  %boat 
			  \begin{scope}
			  		% Parabola
			  		\coordinate (start) at (4, 1.125);
			  		\coordinate (end) at (5, 0.72);
			  		\coordinate (control) at (4.5, 0.78);
			  		
			  		\draw[thick,line width=1.5pt] (start) .. controls (control) .. (end);

			  		% Parabola
			  		\coordinate (start) at (4, 1.125);
			  		\coordinate (end) at (5.08, 1.36);
			  		\coordinate (control) at (4.5, 1.35);
			  		
			  		\draw[thick,line width=1.5pt] (start) .. controls (control) .. (end);
			  	
			  		% Line segment
			  		\draw[thick,line width=1.5pt] (5.025, 0.69) -- (5.065, 1.38);
			  \end{scope}
			  
			 \end{tikzpicture}
		
		 % Tack 2
		\begin{tikzpicture}[scale=0.9]
			
			\coordinate (A) at (4.553, 1.224);
			\coordinate (B) at (5.447,0.776);
			
			\coordinate (C) at ($(A)!.5!(B)$);

			\draw (A) -- (B);
			\pgfmathanglebetweenpoints{\pgfpointanchor{C}{center}}{\pgfpointanchor{A}{center}}
			\let\StartAngle\pgfmathresult
			\pgfmathanglebetweenpoints{\pgfpointanchor{C}{center}}{\pgfpointanchor{B}{center}}
			\let\EndAngle\pgfmathresult

			\draw[fill = red!40] (A) let \p1 = ($(C) - (A)$),
			\n1 = {-veclen(\x1, \y1)}
			in 
			arc [start angle= \EndAngle , end angle= \StartAngle , radius=\n1];
			
			\draw[dashed,line width = 0.3 mm] (3.5,5.5) -- (5,1);
			\draw (3.5,5.5) node[draw,circle,fill=blue!20,inner sep=5pt]{$\tgt$};
			\node[] at (1.9,6.3) {\large \textcolor{red}{\textbf{Tack $q=2$}}};
			
			\draw [decorate, line width=.5mm,
			decoration = {calligraphic brace, raise = 20pt}] (3.5,5.5) -- (5,1)
			node[pos=0.45, right=30pt,black]{$r$};
			
			\tkzDefPoints{5/1/A,3.8/4.6/B,2/4/C,1/3/D,6/0.5/E}
			\tkzDrawSegments[vector style, violet!50, line width = 0.5 mm](A,B)
			\tkzDrawSegments[vector style,violet, line width = 0.5 mm](A,C)
			\tkzDrawSegment[vector style,violet!50, line width = 0.5 mm](B,C)
			\tkzDrawSegment[vector style, green!40!gray,line width = 0.5 mm](D,E)
			\tkzMarkAngle[<-, red, size = 1.5,line width = 0.4 mm](C,A,D)
			\tkzLabelAngle[red, pos=1.2](C,A,D){$u$}
			\tkzMarkAngle[->, green!40!gray, size = 2,line width = 0.4 mm](B,A,D)
			\tkzLabelAngle[green!40!gray, pos=1.7, left = 12pt](B,A,D){$\theta$}
			\tkzLabelSegment[violet!40, pos = 0.6, right = 2 pt ](A,B){$dr$}
			\tkzLabelSegment[violet,swap, pos = 0.9, left=4pt](A,C){$\vec{v}$}
			\tkzLabelSegment[violet!50](C,B){$d \theta$}
			\tkzLabelSegment[green!40!gray, pos = 0.15, below left =0.5 pt](D,E){\textbf{Wind}}
			 %boat 
			\begin{scope}
				% Parabola
				\coordinate (start) at (4.3, 1.7);
				\coordinate (end) at (4.8, 0.76);
				\coordinate (control) at (4.43, 1.2);
				
				\draw[thick,line width=1.5pt] (start) .. controls (control) .. (end);

				% Parabola
				\coordinate (start) at (4.3, 1.7);
				\coordinate (end) at (5.235, 1.217);
				\coordinate (control) at (4.82, 1.58);
				
				\draw[thick,line width=1.5pt] (start) .. controls (control) .. (end);
			
				% Line segment
				\draw[thick,line width=1.5pt] (4.78, 0.78) -- (5.234, 1.235);
			\end{scope}
	\end{tikzpicture}}
	\caption{System dynamics relative to the wind and relative to the target $\tgt$.}
	\label{fig:speed_and_tacks}
\end{figure*}

% COMMENT REMOVED
Following \cite{Ferretti2019, Miles2021}, we assume the strength of the wind is fixed but its direction 
(measured counterclockwise from the $y$-axis)
undergoes a Brownian drift/diffusion process:
\begin{equation}\label{eq:Wind_evo}
	\difd \phi = a \difd t + \sigma \difd B,
\end{equation}
where $\phi(t)$ denotes the current upwind direction, $a$ is a constant drift, $\sigma$ is the diffusion coefficient, and $B$ is a standard Brownian motion.
% COMMENT REMOVED
The state of the system can be represented as $(x,y, q,\phi)$, where $x$ and $y$ encode the boat's current position, while $q \in \{1,2\}$ is its current ``tack'',
which determines the range of available steering directions.    
Our continuous control is the steering angle, % COMMENT REMOVED
$u \in [0, \pi]$, measured relative to the wind.
In the starboard tack ($q=1$), $u$ is measured counterclockwise from the upwind direction, while in the port tack ($q=2$) it is measured clockwise; so,
the boat's direction of motion relative to upwind is $(-1)^q u.$  The boat's angle-dependent speed $f(u)$ is encoded in the speed profile (often called the ``polar''), which is determined by the geometry of each specific boat.  Fig.~\ref{fig:speed_and_tacks}(a) shows a typical polar used in all numerical tests here and in \cite{Miles2021}.
With these assumptions, the boat's dynamics is described by
\begin{align}
\label{eq:full_BoatDyn1}
    dx &= -f(u) \sin(\phi - (-1)^qu) dt, \\
\label{eq:full_BoatDyn2}
    dy &= f(u) \cos(\phi - (-1)^qu) dt. %\\
% COMMENT REMOVED
\end{align}
For the sake of computational efficiency,
we adopt a dimensional reduction of \eqref{eq:Wind_evo}-\eqref{eq:full_BoatDyn2} suitable when aiming for a circular target $\tgt$ in a domain with no obstacles \cite{Miles2021}.
Assuming a target
at the origin, $r = \sqrt{x^2 + y^2}$ represents the boat's distance to the center of $\tgt,$ while $\theta = \phi + \mathrm{atan2}(-x,-y)$ encodes 
the upwind direction measured counterclockwise from the line connecting the boat to the center of $\tgt;$  see Fig.~\ref{fig:speed_and_tacks}(b).
% COMMENT REMOVED
This results in the system dynamics:
\begin{align}\label{eq:Dynamics_r_theta}
	\difd r &= -f(u) \cos(\theta - (-1)^q u) \; \difd t, \notag \\ 
	\difd \theta &= \Bigg[\dfrac{f(u)}{r} \sin(\theta - (-1)^q u) + a \Bigg]\difd t + \sigma  \difd B.
\end{align}
% COMMENT REMOVED
Whenever we focus on the \textit{deterministic portion} of % COMMENT REMOVED
dynamics, we will use 
$r_\difd(\theta,q,u)$ and $\theta_\difd(r,\theta,q,u;a)$
% COMMENT REMOVED
to denote the respective factors multiplying $\difd t$ in the above equations.
% COMMENT REMOVED

% !TEX root = Main_text.tex

\section{Stochastic Optimal Control}\label{sec:Sto_OC}

Let $\bxi(t) = (r(t),\theta(t))$ denote the continuous component of system state at the time $t$. We define $\Omega := [0,\Rmax] \times [0,2\pi) \times \{1,2\}$ as the 
full state space (with the last component encoding the tack $q$).
We will further use $\tilde{q} = 3 - q$ to refer to ``the opposite tack'' and $\Upsilon := [0,\pi] \cup \{\blacktriangle\}$ to define the policy domain (the prescribed steering angle or the special tack-switch action $\blacktriangle$). 

Given a starting configuration of the boat $(\bxi(0) = \hat{\bxi}, q(0) = \hat{q})$ and a feedback control policy $\mu: \Omega \to \Upsilon,$ the main quantity of interest is the random time to target $T_{\mu}(\hat{\bxi},\hat{q}) = \inf \{ t>0 ~|~ \bxi(t) \in \tgt, \; \bxi(0) = \hat{\bxi}, \; q(0) = \hat{q}\}$  
% COMMENT REMOVED
 based on that $\mu.$  Due to the hybrid nature mentioned in section \ref{sec:Sys_Dynamics}, $T_{\mu}(\hat{\bxi},\hat{q})$ is the sum of 
the time spent in both steering and tack-switching.
A typical {\em risk-neutral} approach is to define the value function as
$v(\hat{\bxi},\hat{q}) = \inf_\mu \E[T_\mu( \hat{\bxi},\hat{q})],$ which can be then recovered by solving a system of quasi-variational HJB-type inequalities.
Assuming each tack-switch takes a fixed time $C$ and the boat stays in place ($f(u) =0$) while switching, as shown in \cite{Miles2021}
this leads to % COMMENT REMOVED
\begin{equation}\label{eq:HJB-risk-neutral}
	\max \left\lbrace  H(r,\theta, q, \nabla v) - \dfrac{\sigma^2}{2} \frac{\partial^2 v}{\partial \theta^2},  \; v - \mathcal{N} v - C \right\rbrace = 0, 
\end{equation}
where $\nabla v = (\frac{\partial v}{\partial r}, \frac{\partial v}{\partial \theta})$ and 
the Hamiltonian is
\begin{equation}\label{eq:Hamiltonian-v}
	H(r,\theta, q,\mathbf{p}) = \max_u \left(-	r_{\mathrm{d}}(\theta,q,u) p_1  - \theta_{\mathrm{d}}(r,\theta,q,u;a)p_2 \right) -1.
\end{equation}
In  \eqref{eq:HJB-risk-neutral}, we take the maximum over two alternative courses of action.
The first clause corresponds to the system states, from which 
it is optimal to continue in the current tack and the maximization in \eqref{eq:Hamiltonian-v} selects the optimal steering angle.
On the other hand, in the second clause $\mathcal{N} v$ encodes the expected remaining time-to-target if we switch to the opposite tack;
i.e., $\mathcal{N} v(\hat{r}, \hat{\theta}, q) = \E[v(\hat{r}, \theta(C) , \tilde{q}) \mid \theta(0) = \hat{\theta}, f=0].$
If we define
$\psi^{}_{r, q} (z) = v(r, z, \tilde{q})$ and 
\begin{equation}\label{eq:gaussian-operator}
\mathcal{G}_{\theta} \left[ \psi \right] = \frac{1}{\sigma \sqrt{2 \pi C}} \int_{-\infty}^{\infty} 
% COMMENT REMOVED
\mathrm{e}^{-\mathlarger{\frac{(z - \theta - aC)^2 }{ 2C\sigma^2}}} \psi(z) \, \difd z,
\end{equation}
then this switching operator can be conveniently evaluated as $\mathcal{N} v(r, \theta, q) = \mathcal{G}_{\theta} \left[ \psi^{}_{r, q} \right].$
Since the part of $\Omega$ on which it is better to switch tacks is a priori unknown, this is a problem with a {\em free boundary}.
The optimal feedback policy $\mu_*$ (found by solving \eqref{eq:HJB-risk-neutral} with the boundary condition $v = 0$ on $\tgt \times \{1,2\}$)
captures both the optimal switching states and the optimal steering angles \cite{Miles2021}.
% COMMENT REMOVED

% COMMENT REMOVED
Despite its frequent use, the risk-neutral planning has a significant drawback: it is indifferent to the level of variability in the distribution of times to target. The resulting $\mu_*$ might be impractical if the risk of significantly exceeding 
% COMMENT REMOVED
$\E\left[T_{\mu_*}\right]$
is high (e.g., in right-heavy-tailed distributions).
% COMMENT REMOVED
To address this, we change the perspective and search for a policy $\alpha$ that maximizes the probability of reaching the target before a specified deadline $\shat$. 
We refer to such $\alpha$ as a \textit{risk-aware}  (or, more precisely, as an {\em $\shat$-threshold-aware}) policy.
% COMMENT REMOVED
 
 Letting $\Omega_{\sbar} = \Omega \times [0,\sbar],$ we define our new value function 
% COMMENT REMOVED
\begin{equation}\label{eq:value-fn-thres-aware}
	w(\hat{\bxi}, \hat{q}, \hat{s}) = \sup_{\alpha} \pr \left( T_{\alpha}(\hat{\bxi},\hat{q})  \le \hat{s} \right), \quad 0 \le \hat{s} \le \sbar.
\end{equation}
The supremum is taken over all measurable \textit{threshold-aware} feedback policies $\alpha:\Omega_{\sbar} \to \Upsilon, \; (r,\theta,q,s) \to \alpha(r,\theta,q,s)$.
% COMMENT REMOVED
If we treat $s(0) = \hat{s}$ as an \textit{initial budget}, then the remaining time budget $s(t)$ is strictly decreasing along the path-to-target as $t$ increases ($\dot{s}(t) = -1$ while continuously steering and a negative jump of $C$ units of time when switching tack).
% COMMENT REMOVED

Consequently, the Stochastic Dynamic Programming Principle \cite{Fleming2006} yields an $s$-dependent quasi-variational inequality
\begin{equation}\label{eq:SDPP}
	w(\hat{\bxi}, \hat{q}, \hat{s}) = \max\bigg\{ \sup_{u} \Eu w(\hat{\bxi}, \hat{q}, \hat{s}), \Ec w(\hat{\bxi}, \hat{q}, \hat{s})\bigg\} + o(\tau),
% COMMENT REMOVED
\end{equation}
where 
\begin{align}
	&\Eu \; w(\hat{\bxi}, \hat{q}, \hat{s}) = \E_{u,\tau}[w(\bxi(\tau) , \hat{q}, \hat{s} - \tau) \mid \hat{\bxi},\hat{s}], \label{eq:op1}\\
	&\Ec \; w(\hat{\bxi}, \hat{q}, \hat{s}) = \E_{{}_C}[w(\bxi(C) , \tilde{q}, \hat{s} - C) \mid \hat{\bxi},\hat{s}, f=0]. \label{eq:op2}
\end{align}
% COMMENT REMOVED
Similarly to the structure of \eqref{eq:HJB-risk-neutral}, 
$\Eu w$ refers to
% COMMENT REMOVED
the best probability of reaching the target before the deadline if we stay on the current tack for a small time $\tau$ while 
% COMMENT REMOVED
$\Ec w$ is the best probability if we immediately switch tacks. 

From the stochastic Taylor expansion of Eqn.~\eqref{eq:SDPP}, one can show that, if $w(r,\theta,q,s)$ is sufficiently smooth, it must satisfy
\begin{align}\label{eq:HJB-risk-aware}
	&\max \left\{ \max_{u \in [0,\pi]} \left(  \nabla w^\top \bxi_{\difd}+ \dfrac{\sigma^2}{2} \frac{\partial^2 w  }{\partial \theta^2} -\frac{\partial w}{\partial s}\right)   , \Ec w - w \right\} = 0,
\end{align}
where $\nabla w = (\partial w/\partial r, \; \partial w/\partial \theta)$ and $\bxi_{\mathrm{d}}(r,\theta, q,u) = \Big(r_{\mathrm{d}}(\theta,q,u), \; \theta_{\mathrm{d}}(r,\theta,q,u;a)\Big)$.
% COMMENT REMOVED

As in the risk-neutral case, if we define $\psi_{r,q,s}(z) = w(r,z,\tilde{q},s-C)$, then $\Ec w(r,\theta, q,s) = \mathcal{G}_{\theta}[\psi_{r,q,s}]$ following the definition of operator \eqref{eq:gaussian-operator}.
% COMMENT REMOVED
Note that, in general, the value function does not have to be be smooth or even continuous.
% COMMENT REMOVED
Nonetheless, % COMMENT REMOVED
even a discontinuous value function 
can be often interpreted as a  unique \textit{discontinuous viscosity solution} of a HJB PDE \cite[Chapter 5]{Bardi_book},
which allows recovering the optimal feedback policy $\alpha_*(r,\theta,q,s)$ as 
an argmax of \eqref{eq:HJB-risk-aware}.
% COMMENT REMOVED

% !TEX root = Main_text.tex

\section{Numerical Implementation}\label{sec:Num_Implementaion}
\subsection{Semi-Lagrangian Discretization}\label{sec:SL}
% COMMENT REMOVED
We approximate the solution to \eqref{eq:HJB-risk-aware} for both tacks simultaneously using a semi-Lagrangian discretization \cite{Falcone_book} on a uniform rectangular grid over the $(r,\theta,s)$ space. I.e., $(r_i,\theta_j,s_k) = (i\Delta r, j\Delta \theta, k\Delta s),$ 
where $\Delta r = R_{max} / N_r, \, \Delta r = 2 \pi / N_{\theta}, \, \Delta s = \bar{s} / N_s,$
% COMMENT REMOVED
while $i = 0, ..., N_r, \, j = 1, ..., N_{\theta}$ (due to periodic boundary conditions), and $k = 0, ..., N_s$. 
We will use $W_{i,j}^{k,q} \approx w(r_i,\theta_j,q,s_k)$ to denote the discretized approximate solution at $(r_i, \theta_j,q,s_k)$.
% COMMENT REMOVED
Recall from 
% COMMENT REMOVED
section \ref{sec:Sto_OC} that $s(t)$ is strictly decreasing along the path-to-target. We can thus \textit{causally} march from smaller $s$ to larger $s$ and compute the value function from $s=0$ to $s=\sbar$ in a single sweep. In particular, we choose $\tau = \Delta s$ when solving Eqn.~\eqref{eq:op1} so that we can march from the $s_{k-1}$-slice to the $s_{k}$-slice. 
The expectations in both \eqref{eq:op1} and \eqref{eq:op2} can be approximated using Gauss-Hermite quadratures (GHQ), but the details are somewhat different due to the contrast in 
elapsed times.
% COMMENT REMOVED

To approximate $\Eu,$ we use a first-order weak approximation \cite{Ferretti2010,Kloeden1992} of the distribution of the Brownian increment $\Delta B_\tau$ for $\tau$ time units.   
% COMMENT REMOVED
Starting from any gridpoint $(r_i,\theta_j,s_k)$ and
using any admissible steering angle $u$, we % COMMENT REMOVED
consider two
possible locations of $\bxi(\tau, u)$ in the $s_{k-1}$-slice: 
$$\bxi_{i,j,u}^{\pm} \; = \; (r_i + \tau r_\difd(u), \, \theta_j + \tau \theta_\difd(u)  \pm \sigma\sqrt{\tau} ).$$
Averaging the value function at these points is equivalent to a two-node GHQ approximation of \eqref{eq:op1}; i.e.,
\begin{equation}\label{eq:discrete-op1}
	\Mu W_{i,j}^{k, q} \; = \; \dfrac12 \Big( W^{k-1,q} \left( \bxi_{i,j,u}^{+} \right) + W^{k-1,q} \left( \bxi_{i,j,u}^{-} \right) \Big) .
\end{equation}
Since these $\bxi_{i,j,u}^{\pm}$ are usually not gridpoints, 
we implement \eqref{eq:discrete-op1} by using a \textit{bi-cubic ENO}  interpolation \cite{Shu1998} with a $2\pi$-periodicity in $\theta$. 
% COMMENT REMOVED
We adopt a two stage process for finding the optimal $u_*$ that maximizes $\Mu$ : 
first, we perform a direct comparison over a grid of angle values $\mathcal{U}$ and identify an interval containing the best  $u_{\text{\tiny $\bigstar$}} \in \mathcal{U}$; we then perform a Golden Section Search (GSS) over that interval to obtain a more accurate approximation of $u_*.$

The accuracy of \eqref{eq:discrete-op1} improves under grid refinement\footnote{Under mild technical conditions,
semi-Lagrangian schemes have been proven to converge
to the discontinuous viscosity solution of first-order HJB PDEs on every compact set away from discontinuity 
\cite{Bardi1994fully, Bardi1999numerical}.
For the second-order HJBs, the method closest to ours has been studied (with rigorous error estimates) in \cite{assellaou2015error} 
but without hybrid dynamics or degenerate parabolicity.
While our setting is more general, the numerical results in Section \ref{sec:Num_Experiments} and online repository 
provide strong evidence of convergence.  A rigorous proof of numerical convergence to the discontinuous viscosity solution of PDE \eqref{eq:HJB-risk-aware} remains an open problem to be addressed in the future.
% COMMENT REMOVED
} 
since the diffusion time $\tau = \Delta s \rightarrow 0.$
Finding a good approximation for  $\Ec$  is a bit harder since the diffusion time $C$ is constant.  To address this, one can use a higher order accurate GHQ; e.g., our implementation uses  a version with 3 GH nodes  
% COMMENT REMOVED
\begin{equation}\label{eq:eta}
	\eta_{j,m} \; = \; \theta_j + aC + \sigma\sqrt{2C}x_m, \quad m\in \{1,2,3\}, 
\end{equation} 
where $x_m$ are the roots to the 3rd Hermite polynomial.  Assuming that $s_k \geq C$, we use 
\begin{equation}\label{eq:discrete-op2}
	\Mc W_{i,j}^{k,q} \; = \; \dfrac{1}{\sqrt{\pi}} \sum_{m=1}^{3} \gamma_m \, W \left( r_i, \eta_{j,m}, \tilde{q}, s_k- C \right),
\end{equation} 
where $\gamma_m$'s are the weights of the third GHQ.  We choose $\Delta s$ to be a fraction of $C$, ensuring that $s_k- C = s_l$ for some $l < k.$
But $\eta_{j,m}$ are usually not multiples of $\Delta \theta$ and we use a 1D periodic cubic ENO interpolation to evaluate \eqref{eq:discrete-op2}.

The grid value is then computed as  $$W_{i,j}^{k,q} = \max \left(  M_{u_*, \tau} W_{i,j}^{k, q}, \,   \Mc W_{i,j}^{k,q} \right),$$ and we
recover the optimal % COMMENT REMOVED 
steering/switching policy $\alpha_*(r_i, \theta_j, q, s_k)$   as a by-product.
Our full method is summarized in Algorithm \ref{algo:1} using the target radius $\Rtgt,$  
 the maximum sailboat speed $f_{\mathrm{max}},$ and\\ $\Xi = \{(i\Delta r, j\Delta \theta) \mid i = 0,\ldots,N_r,\; j = 0,\ldots, N_{\theta}\}$, 

% COMMENT REMOVED

% COMMENT REMOVED
\begin{algorithm2e}\label{algo:1}
	\SetAlgoLined
	\For {$s_k = k\Delta s, \; k = 0,1,\ldots N_s$}{
    	\For {$\text{every }\bxi_{i,j} \in \Xi % COMMENT REMOVED
		\; \text{ and } \; q \in \{1,2\}$}{
			\eIf{$(r_i - \Rtgt)/f_{\mathrm{max}} > s_k$}{
               	$W_{i,j}^{k,q} \leftarrow 0$\;   % COMMENT REMOVED
            }{
% COMMENT REMOVED
					$W_{i,j}^{k,q} \leftarrow \max_{u} M_{u,\tau}$\;
					\If {$s \geq C$}{
						$W_{i,j}^{k,q} \leftarrow \max \big( W_{i,j}^{k,q}, \, \Mc W_{i,j}^{k,q} \big)$\;
					}					                    	
                    }
                
               }
           }
   \caption{Risk-aware value function computation}
\end{algorithm2e}

\subsection {Trajectory synthesis and ECDF generation}\label{sec:trajectory_syn}

The above PDE solution process yields the optimal threshold-aware policy in feedback form, with the optimal action $\alpha_*(r_i,\theta_j,q,s_k)$ stored at each gridpoint in $\Xi$ and for a range of deadlines ($k = 0, ..., N_s$).  
To recover a sample path-to-target from any specific initial configuration $(\hat{r}, \hat{\theta}, \hat{q})$ and the intended deadline $\hat{s}$, we use Euler-Maruyama scheme \cite{Kloeden1992} with a fixed time step $\Delta t$ on Eqns.~\eqref{eq:Dynamics_r_theta} .
At each time step, we normally use the optimal steering/switching action from the policy recorded for the nearest gridpoint. % COMMENT REMOVED
% COMMENT REMOVED
But the threshold-aware control formulation leaves two ambiguities that have to be resolved in the implementation.
First, if $w(\hat{r},\hat{\theta},\hat{q}, \hat{s}) = 1,$ the current $\hat{s}$ may be more than sufficient to reach the target with probability one
and the actions taken until the remaining time budget $s$ becomes ``barely sufficient'' are not important.  To address this, we use a ``\textit{Deadline-Upgrade}'' approach,
decreasing the initial time-budget to $\hat{s} = \min \left\{s_k \mid w(\hat{r},\hat{\theta},\hat{q}, s_k) = 1 \right\}.$
Second, if during a simulation the sailor is unlucky and later finds herself with $w(\hat{r},\hat{\theta},\hat{q}, \hat{s}) = 0,$ the PDE provides 
no guidance on what to do after that (since she will now definitely miss the original deadline).  
Rather than dismiss such simulations as complete failure, from there on we simply apply the risk-neutral policy $\mu_*$ recovered from Eqn.~\eqref{eq:HJB-risk-neutral}.

Empirical cumulative distribution function (ECDF) for both $\alpha_*$ and $\mu_*$ are obtained through Monte Carlo simulations,
with sample paths generated % COMMENT REMOVED
starting from a specific $(\hat{r},\hat{\theta},\hat{q}, \hat{s})$
under different realizations of wind evolution. 
% COMMENT REMOVED
The ECDFs for the total time-to-target are obtained through the the Kaplan-Meier estimate \cite{Kaplan1958} using the \texttt{MATLAB}'s built-in function \texttt{ecdf()}.
% COMMENT REMOVED

% !TEX root = Main_text.tex
\section{Numerical Experiments}\label{sec:Num_Experiments}
\begin{figure*}[ht]
	\centering
	\subfigure[][]{%
		\label{fig:risk-aware-policies_a}%
		\includegraphics[scale=0.47]{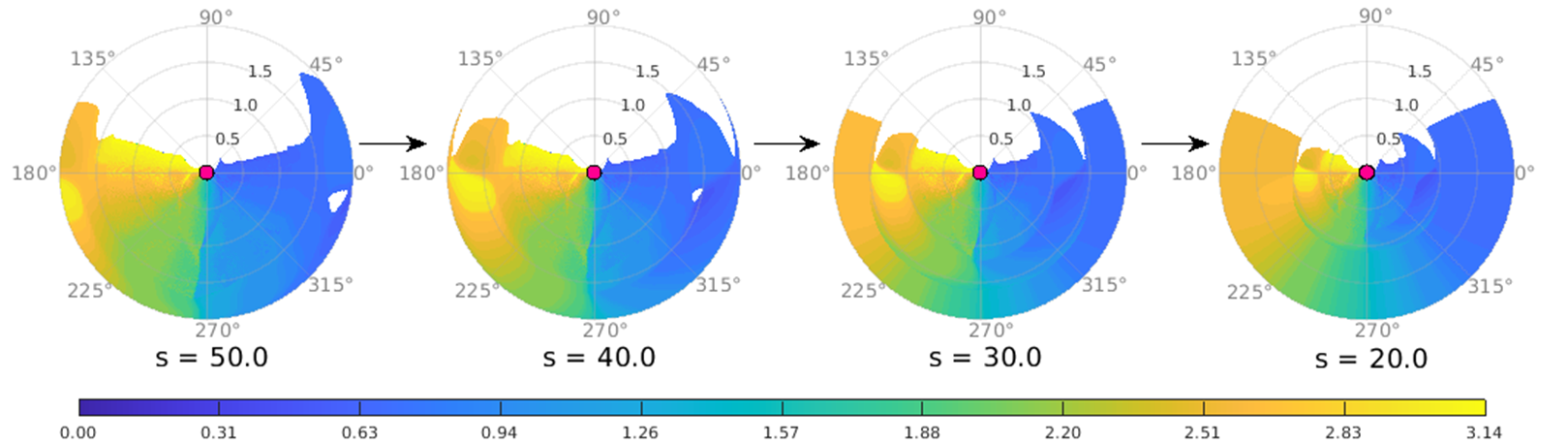}}
	\subfigure[][]{%
		\label{fig:risk-neutral-policies}%
		\includegraphics[scale=0.47]{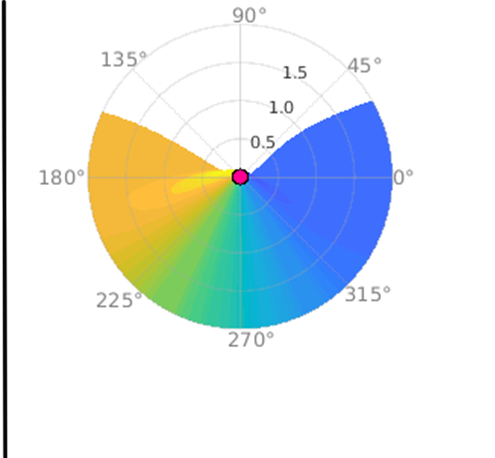}}
	\subfigure[][]{
		\label{fig:risk-aware-policies_b}%
		\includegraphics[scale=0.47]{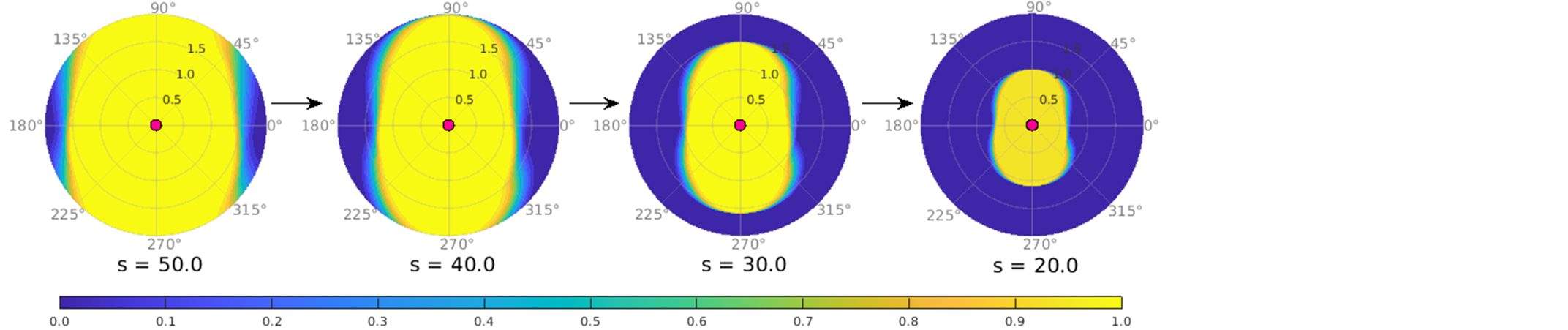}}
	
	\caption[A set of three subfigures.]{\textbf{Representative $s$-slices of risk-aware policy, their corresponding optimal probability of reaching the target $\tgt,$ and the risk-neutral policy}: (a) risk-aware optimal policy $\alpha_*$; (b) risk-neutral optimal policy $\mu_*;$ (c) optimal probability of reaching $\tgt$  associated with (a). All shown for the starboard tack $q=1$ only and in relative $(r,\theta)$ coordinates. 
		In all figures, the target $\tgt$ is shown as a magenta disk in the center.
		In (a) and (b), it is optimal to switch to $q=2$ from wherever the space is left \textit{blank}. 
		Otherwise, it is optimal to stay with $q=1$ and the best steering angle is shown in color (with the same colorbar used in both subfigures). 
% COMMENT REMOVED
		}
	\label{fig:risk-aware-policies}%
\end{figure*}

\begin{figure*}[ht]
	\centering
	\subfigure[][]{%
		\label{fig:CDF_a0_sig005_a}%
		\includegraphics[scale=0.43]{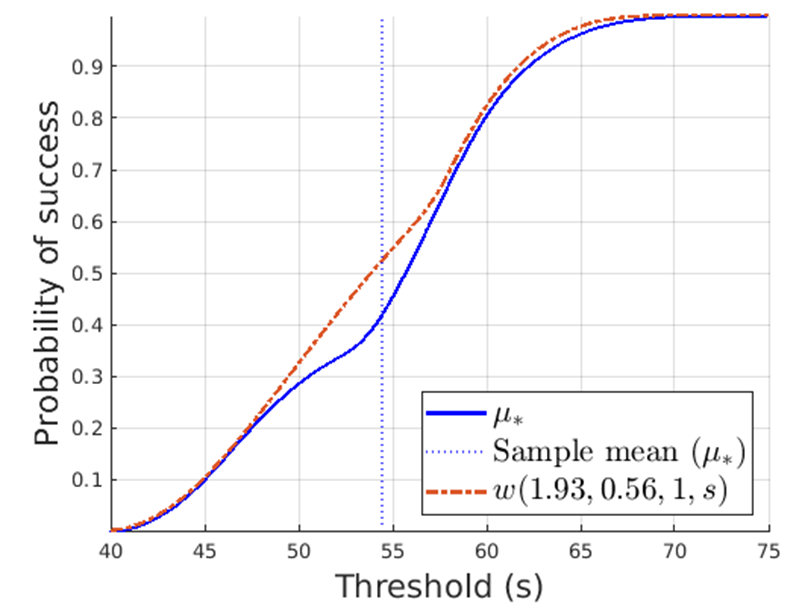}}
	\subfigure[][]{
		\label{fig:CDF_a0_sig005_b}%
		\includegraphics[scale=0.43]{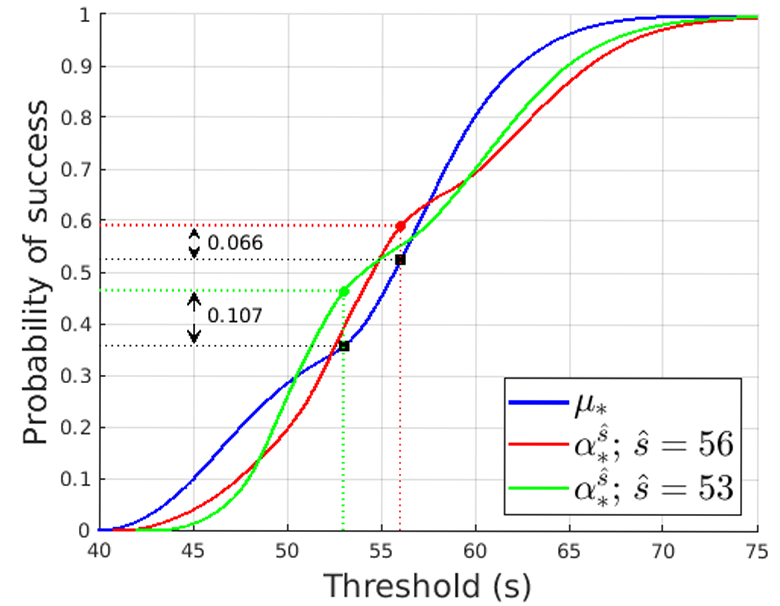}}
	\subfigure[][]{
		\label{fig:CDF_a0_sig005_c}%
		\includegraphics[scale=0.44]{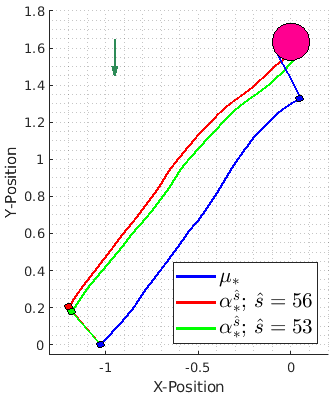}}
	
	\caption[A set of three subfigures.]{\textbf{Sailing against the wind: a comparison between the risk-aware and risk-neutral approaches with $(a=0, \, \sigma = 0.05)$ starting from $(\hat{r},\hat{\theta},\hat{q}) = (1.93, 0.56, 1)$} Subfigures:
	(a) ECDF (empirical cumulative distribution function) generated with the optimal risk-neutral policy $\mu_*$ (\textit{solid blue}) vs. the $s$-dependent risk-aware optimal probability of success 
	$w(\hat{r},\hat{\theta},\hat{q},s)$ (\textit{dash-dotted orange}); 
	(b) ECDFs of the random total time to target generated with different polices; 
	(c) Sample sailboat trajectories in the \textit{absolute $xy$-coordinates} generated with different polices under the same random wind path.
	The target set is plotted as a magenta disk at the top, and top-left \textit{dark green arrow} encodes the initial wind direction $\phi(0) = 0.$ 
	Trajectory colors correspond to the polices used to generate ECDFs in (b). 
	The colored dots indicate the tack-switching points for respective trajectories.
	Observed sample means for the arrival time $T$ are  $54.46, \, 55.57,$ and $55.95$ for the policies $\mu_*, \, \alpha_*^{53},$ and $\alpha_*^{56}$ respectively.
	Arrival times for the specific trajectories in (c) are $56.55$ (blue), $53.54$ (green), and $54.07$ (red).
	}
% COMMENT REMOVED
	\label{fig:CDF_a0_sig005}%
\end{figure*}

We use several examples to compare the performance of risk-aware and risk-neutral policies.
In all cases, the value functions are computed on a $1601\times 1601 \times 2$ grid for $(r,\theta,q) \in [0,2] \times [0, 2\pi] \times \{1,2\}.$ 
When solving  \eqref{eq:HJB-risk-aware}, we use $\Delta s = 0.025$ for any preset maximum deadline $\sbar$.
The other parameter values are $\Rtgt = 0.1$, $C=2$, and $f_{\mathrm{max}} = 0.05.$ 
All ECDFs are built through Monte Carlo simulations (see section \ref{sec:trajectory_syn}) with $10^5$ samples and $\Delta t = 0.005$. 

% COMMENT REMOVED

% COMMENT REMOVED
We illustrate our $s$-dependent risk-aware policies for a wind with zero drift ($a=0$) and small diffusivity ($\sigma = 0.05$). 
In Fig.~\ref{fig:risk-aware-policies}(a,c), we show some representative $s$-slices of the risk-aware optimal policies $\alpha_*$ and the corresponding value function $w$ for the starboard tack\footnote{
In the interest of reproducibility, our full source code, additional examples, and movies (for both tacks) will be available from\\ 
\href{https://eikonal-equation.github.io/Threshold-Aware-Sailing-Public}{https://eikonal-equation.github.io/Threshold-Aware-Sailing-Public}.} 
% COMMENT REMOVED
(i.e., the optimal probability of reaching $\tgt$ in less than $s$ units of time if we start with $q=1$ and use $\alpha_*$). 
% COMMENT REMOVED
In Fig.~\ref{fig:risk-aware-policies}(a),
colors indicate the optimal steering angle $u_*$ in the current tack, while the complement 
(left blank) shows all the $(r, \theta)$ configurations at which the immediate tack-switch $\blacktriangle$ is optimal.
We observe that $\alpha_*$ is strongly $s$-dependent and significantly differs from 
the risk-neutral optimal policy $\mu_*$ shown in Fig.~\ref{fig:risk-aware-policies}(b).
The arrows in Fig.~\ref{fig:risk-aware-policies}(a,c) indicate the natural progression when threshold-aware policies are used in practice: 
once we start with a particular deadline $\shat$, our initial time-budget\footnote{
In the following discussion, we use the superscript ($\alpha_*^{\shat}$) to refer to a version of policy $\alpha_*$ implemented with a specific initial time-budget $%s(0)= 
\shat.$
}
$s(0)= \shat$
is progressively decreasing, making it necessary to use $\alpha_*$ from the lower $s$-slices.
This decrease is gradual ($\dot{s}(t) =-1$) while we are steering, but becomes abrupt whenever we decide to tack-switch at some time
$t^{\blacktriangle}$; i.e., $s(t) = s(t^{\blacktriangle}) - C, \quad \forall t \in (t^{\blacktriangle}, t^{\blacktriangle}+C].$

% COMMENT REMOVED

Since $\alpha_*$ is somewhat more complicated to implement in practice, 
it is reasonable to ask whether it is significantly better (in meeting the desired deadlines) than the risk-neutral $\mu_*.$
For any specific starting configuration, the answer can be found by comparing the ECDF of $\mu_*$ with the risk-aware value function  $w$ plotted across the range of $s$ values.
The graph of $w$ will be always above, though often this difference is minimal, making the use of $\mu_*$ a preferred option.
However, in many cases the gap between the two graphs will be more significant for a specific range of $s$ values.
This is illustrated in Fig.~\ref{fig:CDF_a0_sig005_a} for % COMMENT REMOVED
$(\hat{r},\hat{\theta},\hat{q}) = (1.93, 0.56, 1).$
If we are interested in some deadline between $\shat=52$ and $\shat=57,$ the threshold-aware policies provide a noticeable advantage.
For example, our $\alpha_*^{56}$ (red in Fig.~\ref{fig:CDF_a0_sig005_b}) increases $\pr(T \le 56)$ from $52.5\%$ to $58.9\%,$
% COMMENT REMOVED
while
our 
$\alpha_*^{53}$ (green in Fig.~\ref{fig:CDF_a0_sig005_b}) yields a $10.6\%$ improvement in $\pr(T \le 53)$
while increasing $\E[T]$ by less than $2.8\%.$

It is also % COMMENT REMOVED
revealing to examine sample trajectories resulting from each of these policies (shown in Fig.~\ref{fig:CDF_a0_sig005_c} for a particular random realization of wind evolution).  
According to $\mu_*$, the boat starts in the ``wrong'' tack, and thus needs to switch % COMMENT REMOVED
immediately, with another tack-switch (back to $q=1$) almost always needed later % COMMENT REMOVED
to reach $\tgt$.  
This strategy produces the best  $\E[T],$ % COMMENT REMOVED
but makes it hard to reach the target much earlier and does not hedge against the bad outcomes (e.g., $T_{\mu_*} > 58$ in more that 33\% of simulations).
In contrast, the threshold-aware policies make a calculated bet (that the wind direction will soon change to help us), stay with the original $q=1$ at first, and 
% COMMENT REMOVED
reach $\tgt$ with only one tack-switch.  

\begin{figure}[ht]
	\centering
	\subfigure[][]{
		\label{fig:CDF_a005_sig005}
% COMMENT REMOVED
		\includegraphics[scale=0.3]{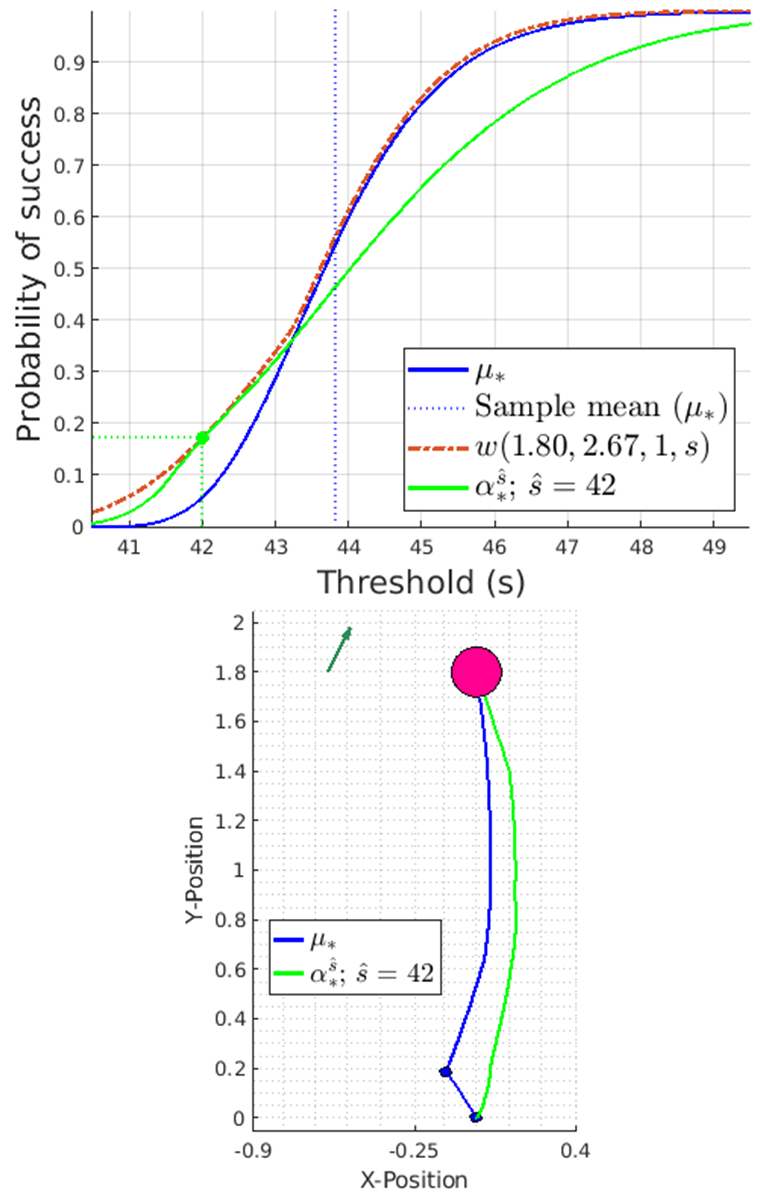}}
	\subfigure[][]{
		\label{fig:CDF_a015_sig005}
% COMMENT REMOVED
		\includegraphics[scale=0.3]{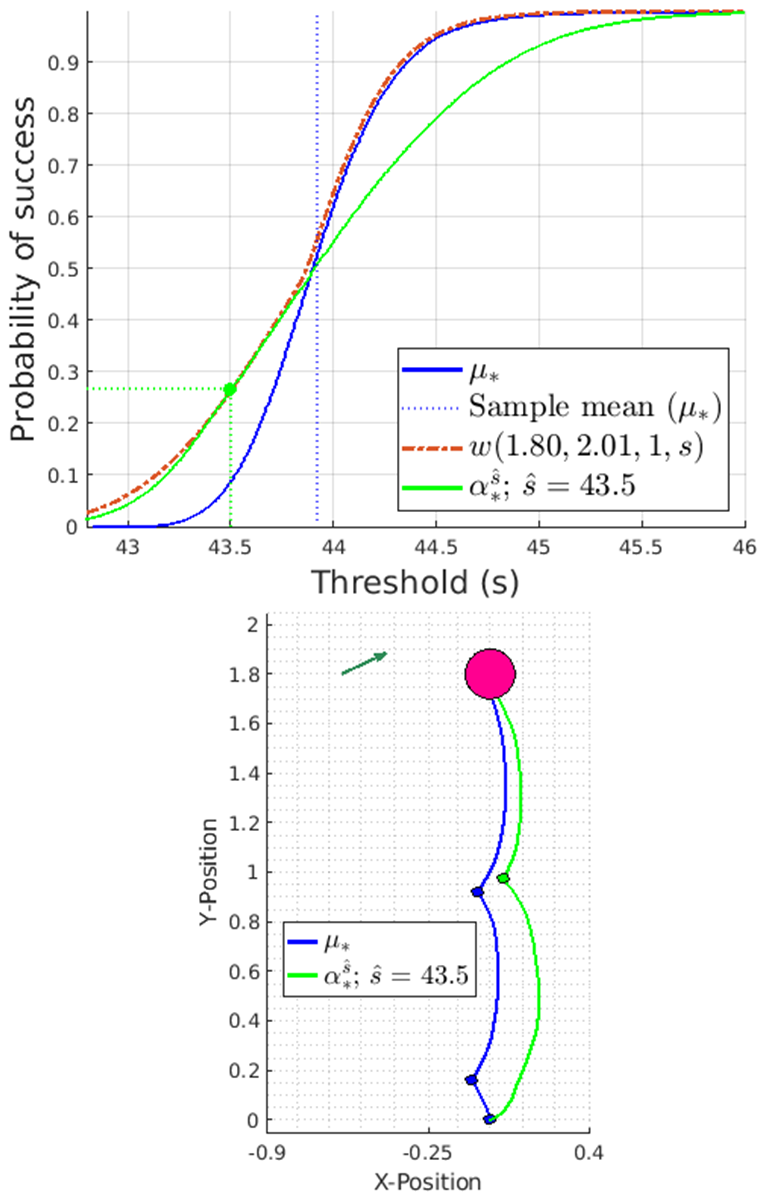}}
			\caption[A set of two subfigures.]{\textbf{Exploiting the wind-drift}: (a) $(a=0.05, \, \sigma = 0.05)$; (b) $(a=0.15, \, \sigma = 0.05)$.
\textbf{Top Row:} ECDF for $\mu_*$ (\textit{solid blue}), the $s$-dependent risk-aware optimal probability of success $w(\hat{r},\hat{\theta},\hat{q})$ (\textit{dash-dotted orange}),
and ECDF for $\alpha_*^{\shat}$ (\textit{solid green}). 
In (a), $(\hat{r},\hat{\theta},\hat{q}) = (1.80,2.67,1)$ and $\shat = 42$. The sample means for $\mu_*$ and $\alpha_*^{\shat}$ are $43.83$ and $44.30.$
In (b), $(\hat{r},\hat{\theta},\hat{q}) = (1.80,2.01,1)$ and $\shat = 43.5$. The sample means for $\mu_*$ and $\alpha_*^{\shat}$ are $43.93$ and $43.97.$
\textbf{Bottom Row}: Two representative % COMMENT REMOVED
paths generated with the same wind evolution (with colors corresponding to respective policies in the top row). 
The \textit{dark green arrow} encodes the initial wind direction. % COMMENT REMOVED
\textit{Time-to-target:} (a) blue: $42.23$, green: $41.44$ ; (b) blue: $43.12$, green: $42.98$. 
% COMMENT REMOVED
In (a) $\mu_*$ led to 2 tack-switches in $99.9 \%$ of simulations, while $\alpha_*$ required none in $99.1\%$ of cases with 2 switches needed in all others.
In (b) $\mu_*$ led to 3 tack-switches in $99.9 \%$ of simulations (with others requiring 4), while $\alpha_*$ required  1 switch in $99.9\%$ of cases with 2 switches needed in the rest.
% COMMENT REMOVED
}
% COMMENT REMOVED
	\label{fig:CDF_different_wind}%
\end{figure}

Larger improvements 
can be % COMMENT REMOVED
similarly realized
with a non-zero wind-drift, particularly when the chosen deadlines are fairly aggressive  (in the left tail of $T_{\mu_*}$ PDF).
In Fig.~\ref{fig:CDF_different_wind}(a) we show such an example with  $(a = 0.05, \, \sigma = 0.05).$  Unlike in Fig.~\ref{fig:CDF_a0_sig005}, here the initial direction of the wind is largely toward the target,
but we are in the wrong initial tack to fully take advantage of this.  The risk-neutral $\mu_*$ prescribes an immediate tack-switch followed by another one a bit later and yields 
a low $\pr(T \le 42) \approx 5.8\%.$  In contrast, the threshold-aware $\alpha_*^{42}$ recognizes, based on the sign of $a,$ that the wind is likely to change in the right direction soon 
and (in the specific wind-evolution example presented in the bottom row of Fig.~\ref{fig:CDF_different_wind}(a)) manages to reach $\tgt$ without any tack-switches at all.
The result of this calculated bet is to almost triple $\pr(T \le 42)$ to $17.2\%$ and make the tack-switches far less common.
Even more dramatic improvements can be obtained when the drift is stronger.
In Fig.~\ref{fig:CDF_a015_sig005} with  $(a = 0.15, \, \sigma = 0.05),$ the threshold-aware policy boosts $\pr(T \le 43.5)$ from $8.8\%$ to $26.6\%,$
largely by reducing the number of tack-switches (in most cases, from 3 switches under $\mu_*$ to only 1 under $\alpha_*^{43.5}$).

% COMMENT REMOVED

% COMMENT REMOVED

We end this 
% COMMENT REMOVED
section
with two caveats.  
First, it is usually impossible to optimize the entire CDF of the arrival time.  As should be clear from Fig.~\ref{fig:CDF_a0_sig005_b}, a policy increasing $\pr(T \le \shat_1)$ might be decreasing 
$\pr(T \le \shat_2)$ even compared to a risk-neutral $\mu_*.$   Typically, each $\alpha_*^{\shat}$ is  only optimal for its particular threshold/deadline $\shat.$  This is why we do not use 
the usual nomenclature of {\em risk-aversion} \cite{WangChapman_2022}  and instead describe or methods as risk (or threshold) {\em aware}.  Second, our decision to revert to $\mu_*$ 
in the ``unlucky'' $\alpha_*$-based simulations (once the time-budget is reduced to zero) is fairly arbitrary and one can certainly use other approaches instead.  However, this choice 
does not affect $\pr\left(T_{\alpha_*^{\shat}} \le \tilde{s} \right)$ for any $\tilde{s} \leq \shat$; thus, the primary goal of threshold-aware policies is still achieved.

% !TEX root = Main_text.tex
\section{Conclusion}\label{sec:Conclusion}
% COMMENT REMOVED
We have introduced a robust (risk/deadline-aware) approach to controlling a sailboat in stochastically evolving wind conditions.
The efficiency of our approach hinges on the numerical method for a pair of quasi-variational HJB-type inequalities, 
which yield deadline-aware policies for all initial configurations and a broad range of deadlines simultaneously.
Numerical experiments demonstrate the advantages of these policies over the traditional risk-neutral approach \cite{Miles2021}, particularly
when it is possible to reduce the number of likely tack-switches.

Several extensions will obviously increase the impact of this approach in the future.  Solving the problem in absolute coordinates will allow for a better modeling of the domain geometry (e.g., accounting for obstacles and other target shapes). Incorporating more realistic wind models and more detailed boat dynamics will be clearly of interest to practitioners.  
Similarly, stochastic differential games might be used to reflect the competitive aspect of sailing races \cite{cacace_stochastic_2019}.  
[E.g., if $T_i$ is a (random) arrival time of the $i$-th competitor, one could try to maximize $\pr \left( T \leq \min_i T_i \right).$]  
% COMMENT REMOVED
In addition, it will be important to explore multi-objective versions
(e.g., Pareto-optimal tradeoffs between $\E[T]$ and $\pr\left(T \leq \shat \right)$) and compare our approach with % COMMENT REMOVED
risk-averse methods that minimize the ``Conditional Value at Risk'' \cite{WangChapman_2022}.

Finally, we hope that a similar threshold-aware approach will prove to be useful in many indefinite-horizon hybrid control applications unrelated to sailing.   
% COMMENT REMOVED

\addtolength{\textheight}{-12cm}   % This command serves to balance the column lengths

%%%%%%%%%%%%%%%%%%%%%%%%%%%%%%%%%%%%%%%%%%%%%%%%%%%%%%%%%%%%%%%%%%%%%%%%%%%%%%%%

%%%%%%%%%%%%%%%%%%%%%%%%%%%%%%%%%%%%%%%%%%%%%%%%%%%%%%%%%%%%%%%%%%%%%%%%%%%%%%%%

% COMMENT REMOVED

%%%%%%%%%%%%%%%%%%%%%%%%%%%%%%%%%%%%%%%%%%%%%%%%%%%%%%%%%%%%%%%%%%%%%%%%%%%%%%%%
\bibliographystyle{IEEEtran}
\bibliography{IEEEabrv,Thres_aware_Sailing_Bib}
\end{document}